\documentclass{tsp-short}
\usepackage{amsfonts}
\usepackage{amssymb}
\usepackage{newlfont}
\usepackage{amsthm}
\usepackage{euscript}
\usepackage{amsmath}
\usepackage{graphicx}
\usepackage{dsfont}
\newtheorem{Proposition}{Proposition}

\newtheorem{Theorem}{Theorem}

\theoremstyle{definition}

\theoremstyle{remark}

\newtheorem{Remark}{Remark}

\begin{document}

\title{On simultaneous  hitting of  membranes by two skew Brownian motions}

\author{Olga V. Aryasova}
\address{Institute of Geophysics, National Academy of Sciences of Ukraine,
Palladina pr. 32, 03680, Kiev-142, Ukraine}
\email{oaryasova@mail.ru}

\author{Andrey Yu. Pilipenko}
\footnote{Research partially supported by Ministry of Education
and Science of Ukraine,
    Grant ¹ F26/433-2008, and  National Academy of Sciences of
    Ukraine, Grant ¹ 104-2008.}
\address{Institute of Mathematics,  National Academy of Sciences of
Ukraine, Tereshchenkivska str. 3, 01601, Kiev, Ukraine}
\email{apilip@imath.kiev.ua}

\date{01/04/2009}
 \dedicatory{}


\begin{abstract}
We consider two depending Wiener processes which have membranes  at zero with different permeability coefficients.  Starting from different points, the processes almost surely do not meet at any fixed point except that where membranes are situated. The necessary and sufficient conditions for the  meeting of the processes are found. It is shown that the probability of meeting is equal to zero or one.
\end{abstract}
\maketitle

\section*{Introduction}
Let $\left(w_1(t)), (w_2(t)\right)_{t\geq 0}$ be a two-dimensional Wiener process with the correlation matrix
$$
B=
\begin{pmatrix}
1&\alpha\\
\alpha&1
\end{pmatrix}t,
$$
where $\alpha\in(-1,1)$ is some constant.

Consider the equations
\begin{eqnarray}
x_1(t)&=&x_1(0)+w_1(t)+\varkappa_1L_{x_1}^0(t),\label{x_1}\\
x_2(t)&=&x_2(0)+w_2(t)+\varkappa_2L_{x_2}^0(t)\label{x_2},
\end{eqnarray}
where $\{\varkappa_1,\varkappa_2\}\in[-1,1]$
are constants,
$$
L_{x_i}^0(t)=\lim_{\varepsilon\downarrow 0}\frac{1}{2\varepsilon}\int_0^t \mathds 1_{[-\varepsilon,\varepsilon]}(x_i(s))ds,\ i=1,2,
$$
is a local time of the process $(x_i(t))_{t\geq0}$  at zero. As is
known (cf. \cite{HS}), each of equations (\ref{x_1}), (\ref{x_2})
has a unique solution, which is a skew Brownian motion. Here
$\varkappa_1$ and $\varkappa_2$ can be treated as coefficients of
permeability. If $\varkappa_1=1$, the part of the process
$(x_1(t))_{t\geq0}$ on the positive semi-axis  is a Wiener process
with reflection at $0$; if $\varkappa_1=-1$, then the part of the
process $(x_1(t))_{t\geq0}$ on the negative semi-axis is  a Wiener
process with reflection at  $0$; if $\varkappa_1\in(-1,1)$, then there
is a semipermeable membrane at $0$.

The aim of the paper is to calculate the probability of
simultaneous hitting of the membranes by the processes
$(x_1(t))_{t\geq 0}$  and $(x_2(t))_{t\geq 0}$. This probability
turns out to be determined by the sign of an expression involving
$\varkappa_1,\varkappa_2,\alpha.$ Besides it is equal to zero or one.

The case of  $\alpha=1$  is studed in \cite{B_coal}, \cite{B_local}.
In \cite{B_coal} it is proved that if
$\{\varkappa_1,\varkappa_2\}\in[-1,1]\setminus\{0\}$ then the processes
$(x_1(t))_{t\geq 0}$ and $(x_2(t))_{t\geq 0}$ meet in  finite time
with probability 1. In  \cite{B_local} it is obtained that if
$x_1(0)=x_2(0)=0,$  $0<\varkappa_1<\varkappa_2<1,$ and
$\varkappa_1>\varkappa_2/(1+2\varkappa_2),$ then for each $t_0>0$ there exists
$t>t_0$ such that $x_1(t)=x_2(t).$   The problem of simultaneous
hitting of the sphere by two Brownian motion with normal
reflection on the sphere is treated in \cite{Le Jan circle}
(two-dimensional case) and \cite{Sheu}.

If there are no membranes, i.e. $\varkappa_1=\varkappa_2=0, \alpha\neq 1$,
then the process  $x(t)=\left(x_1(t), x_2(t)\right), t\geq 0,$ is
a two-dimensional Wiener process.   It reaches any fixed point
$x_0\in\Bbb R^2,$ $ x_0\neq x(0),$ with probability 0. In
particular this implies that the process $(x(t))_{t\geq 0}$ almost
surely does not hit any fixed point  except the points at which at
least one membrane is situated.

There is  one more problem of stochastic analysis where the study
of simultaneous membrane visitation arises naturally. Assume that
we are attending to construct a flow generated by stochastic
differential equation with a semipermeable membrane located on a
hyperplane \cite{Portenko}. Note that there is no general results
on existence and uniqueness of a strong solution to such equations
in multidimensional space. In order to construct the flow on  some
probability space it is sufficient to construct a sequence of
consistent (weak) $n$-point motions \cite{Darling}. One-point
motion can be constructed by N.Portenko's methods \cite{Portenko}.
There are no general results on weak uniqueness for two-point
motion, when both points start from the membrane. However, if the
simultaneous visitation the membrane has a probability 0, then
there is a hope to construct $n$-point motion using localization
at the neighborhood of membrane. Unfortunately, the results of the
article show that synchronous hitting the membrane is quite
natural.

\section{Transformation of the processes}

The pair of the processes $(x_1(t), x_2(t))_{t\geq 0}$ can be thought off as a new process in Euclidean space $\Bbb R^2$ with membranes on the straight-lines
$S_1=\{x_2=0\}$ ³  $S_2=\{x_1=0\}$. The membranes act in the normal direction
$\nu_1=(0,1)$ and $\nu_2=(1,0)$ to $S_1$ and $S_2$
respectively.

Let us make a coordinate transformation defined by the linear operator
$$
A=B^{-1/2}=\frac{1}{c}
\begin{pmatrix}
a&b\\
b&a
\end{pmatrix},
$$
where
\begin{eqnarray*}
a&=&\sqrt{1-\alpha}+\sqrt{1+\alpha},\\
b&=&\sqrt{1-\alpha}-\sqrt{1+\alpha},\\
c&=&2\sqrt{1-\alpha^2}.
\end{eqnarray*}
As a result we get a new process $(\tilde x_1(t), \tilde x_2(t))_{t\geq 0}.$
From  (\ref{x_1}), (\ref{x_2}) we see that its trajectories are
solutions of the following equations
\begin{eqnarray}
\tilde x_1(t)&=&\tilde x_1(0)+\tilde w_1(t)+\varkappa_1\frac{a}{c}L_{x_1}^0(t)+\varkappa_2\frac{b}{c}L_{x_2}^0(t),\label{tx_1}\\
\tilde x_2(t)&=&\tilde x_2(0)+\tilde w_2(t)+\varkappa_1\frac{b}{c}L_{x_1}^0(t)+\varkappa_2\frac{a}{c}L_{x_2}^0(t)\label{tx_2},
\end{eqnarray}
where $\tilde w_1(t)={a}/cw_1(t)+b/cw_2(t),\
\tilde w_2(t)={b}/cw_1(t)+a/cw_2(t), \ t\geq 0.$ It is easily seen that
the correlation matrix of the vector $(\tilde w_1(t),\tilde w_2(t))$ is as
follows
$$\begin{pmatrix}
1&0\\
0&1
\end{pmatrix}t.
$$
This yields that    $(\tilde w_1(t))_{t\geq 0}$ and $(\tilde w_2(t))_{t\geq 0}$ are independent Wiener processes.

Denote by $S_1^{\,\prime}$ and $S_2^{\,\prime}$ the images of
$S_1$ and $S_2$ under the transformation defined by the matrix
$A$. Then equations (\ref{tx_1}), (\ref{tx_2}) can be rewritten in
the form
\begin{eqnarray}
\tilde x_1(t)&=&\tilde x_1(0)+\tilde w_1(t)+\varkappa_1\frac{a}{c}L_{\tilde x}^{S_1^{\,\prime}}(t)+\varkappa_2\frac{b}{c}L_{\tilde x}^{S_2^{\,\prime}}(t),\label{tx_1 prime}\\
\tilde x_2(t)&=&\tilde x_2(0)+\tilde w_2(t)+\varkappa_1\frac{b}{c}L_{\tilde x}^{S_1^{\,\prime}}(t)+\varkappa_2\frac{a}{c}L_{\tilde x}^{S_2^{\,\prime}}(t)\label{tx_2prime},
\end{eqnarray}
where $L_{\tilde x}^{S_i^{\,\prime}}$ is a symmetric local time of the process $(\tilde x(t))_{t\geq 0}$ on the straight-line $S_i^{\,\prime}$ that is
\begin{equation}
L_{\tilde x}^{S_i^{\,\prime}}(t)=\lim_{\varepsilon\downarrow0}\frac{1}{2\varepsilon}\int_0^t \label{L^c_k} \mathds 1_{A_{\varepsilon}^i}(\tilde x(s))ds, \  i=1,2,
\end{equation}
$$
A_\varepsilon^i=\{x\in\Bbb R^2:\exists y\in S_i^{\,\prime}, s\in[-1,1] \mbox{ such that } x=y+\varepsilon s \nu_i^{\,\prime}\},
$$
$\nu_i^{\,\prime},\  i=1,2, $ is the image of $\nu_i$ under the transformation defined by the matrix $A$.

\section{On hitting of zero by the  Wiener process on the plane with membranes on rays with a common endpoint}

A Wiener process in $\Bbb R^2$ with membranes on rays $c_1,\dots,\ c_n$ having a common endpoint was investigated in \cite{APi}.
Let $(r,\varphi),r\geq0,\varphi\in[0,2\pi),$ be polar coordinates in
$\mathbb{R}^2$ and let
$$
c_k=\{(r,\varphi):r\geq0,\ \varphi=\varphi_k\},
$$
where $0\leq\varphi_1< \dots<\varphi_n<2\pi$. Put $\varphi_{n+1}=\varphi_n,$  $  \xi_k=\varphi_{k+1}-\varphi_k, k=1,2,\dots,n$, and $\xi_0=\xi_n$.

\begin{figure}[h]
  \includegraphics[width=9 cm]{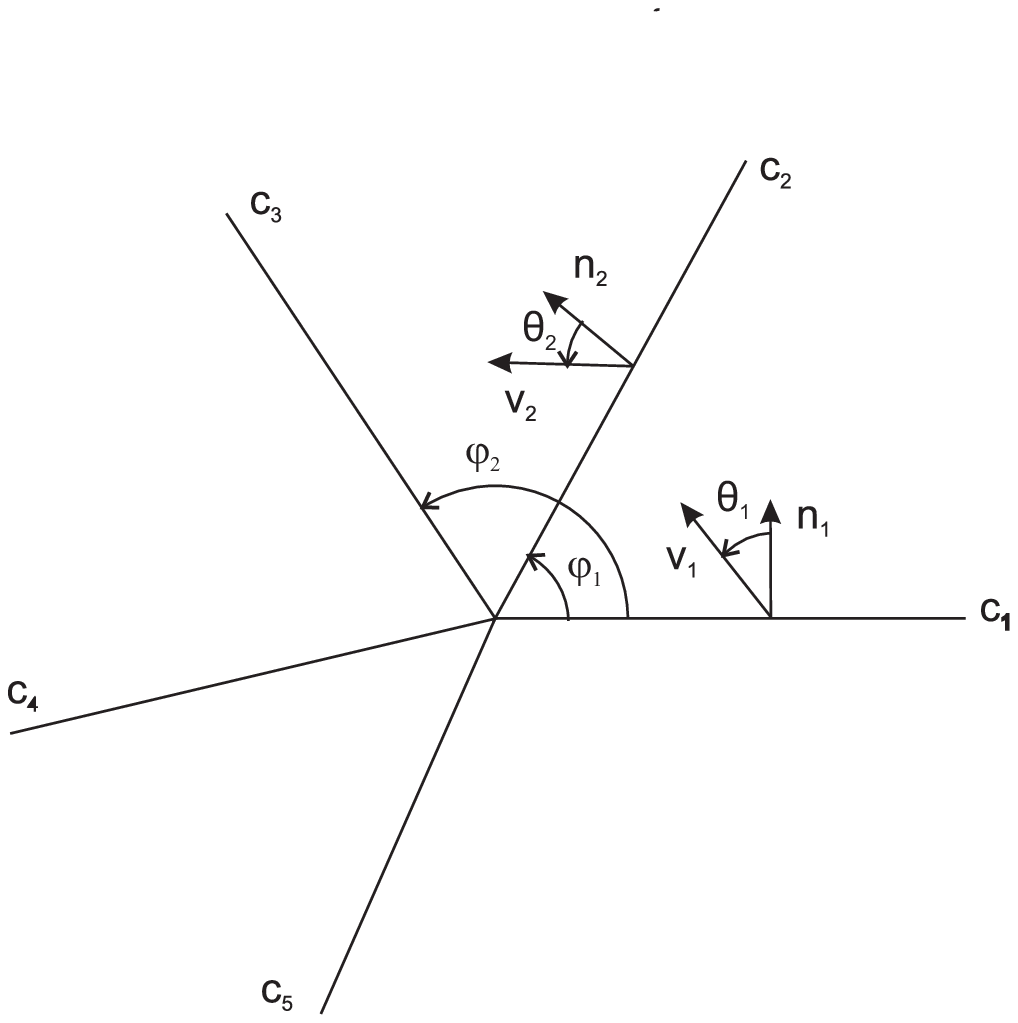}\\
  \caption{}
\end{figure}

Denote by  $n_k, k=1,\dots,n,$ the unit vector normal to $c_k$ that points anticlockwise
and let $v_k$ be a vector in $\mathbb
R^2$ such that $(v_k,n_k)=1$. The angle between $n_k$ and $v_k$
denoted by $\theta_k\in (-\frac{\pi}{2},\frac{\pi}{2})$ is referred to as a positive if and only if $v_k$ points towards the origin.
Let  $\gamma_k, |\gamma_k|\leq1,\ k=1,\dots,n,$ be the membrane permeability coefficients. The case of  $n=5,\theta_1>0,\theta_2>0$ is shown in Fig. 1.

It was proved in  \cite{APi} that there exists a unique strong solution to the equation
\begin{equation}
dx(t)=dw(t)+\sum_{k=1}^n \gamma_k v_k dL_x^{c_k}(t) \label{main}
\end{equation}
in $\Bbb R^2$ with the initial condition
$x(0)=x^0,\ x^0\in\Bbb R^2,$ up to the time  $\zeta,$ where $\zeta=+\infty$ or $x(\zeta-)=0.$

Further on we make use of the following  Proposition on hitting
$0$ or $\infty$ by the process $(\tilde x(t))_{t\geq 0}$ (cf.
\cite{APi}).

\begin{Proposition}
Let $\gamma_k\in[-1,1], \ k=1,\dots,n$, and
let the Markov chain with the state-space $\{1,\dots,n\}$ and the
transition matrix
$$\begin{pmatrix}
0 & \tilde p_1 & 0&0&\hdots&0&0&0&\tilde q_1\\
\tilde q_2 & 0&\tilde p_2 & 0&\hdots &0&0&0&0\\
\hdotsfor[3]{9}\\
0&0&0&0&\hdots&0&\tilde q_{n-1}&0&\tilde p_{n-1}\\
\tilde p_n&0&0&0&\hdots&0&0&\tilde q_n&0
\end{pmatrix},
$$
where
\begin{eqnarray}
\tilde p_k&=&\frac{(1+\gamma_k)\xi_{k-1}}{(\xi_{k-1}+\xi_k)+\gamma_k(\xi_{k-1}-\xi_k)}\label{tp_k}, \\
\tilde q_k&=&\frac{(1-\gamma_k)\xi_{k}}{(\xi_{k-1}+\xi_k)+\gamma_k(\xi_{k-1}-\xi_k)}\label{tq_k},
\end{eqnarray}
has a unique invariant distribution  $(\pi_k)_{k=1}^n$.

Then if $\sum_{k=1}^n \gamma_k\pi_k
\frac{\xi_{k-1}\xi_k}{{(\xi_{k-1}+\xi_k)+\gamma_k(\xi_{k-1}-\xi_k)}}\tan\theta_k>0$
then the process $(\tilde x(t))_{t\geq 0}$  hits the origin almost
surely; \newline if $\sum_{k=1}^n \gamma_k\pi_k
\frac{\xi_{k-1}\xi_k}{{(\xi_{k-1}+\xi_k)+\gamma_k(\xi_{k-1}-\xi_k)}}\tan\theta_k\leq
0$ then the process $(\tilde x(t))_{t\geq 0}$  does not hit the origin
a.s.
\end{Proposition}

\section{The main result}

Let us formulate our problem in  terms of the previous Section.
Let $ c_1, c_2, c_3, c_4$ be the images of the rays
$[0,\infty)\times\{0\},$ $ \{0\}\times[0,\infty)$,
$(-\infty,0]\times\{0\}$, $ \{0\}\times(-\infty,0]$ under the
linear transformation $A$. The images of $\nu_1=(0,1)$ and
$\nu_2=(1,0)$ are the vectors  $\nu_1^{\prime}=(a/c,b/c)$ and
$\nu_2^{\,\prime}=(b/c,a/c)$. Denote by  $\xi $ the angle between
them. Then
$$
\cos\xi=\frac{(\nu_1^{\,\prime},\nu_2^{\,\prime})}{|\nu_1^{\,\prime}|\cdot|\nu_2^{\,\prime}|}=\frac{2ab}{a^2+b^2}=-\alpha.
$$
The case of $\alpha<0$ is shown in  Fig. 2.
\begin{figure}[h]
  \includegraphics[width=7 cm]{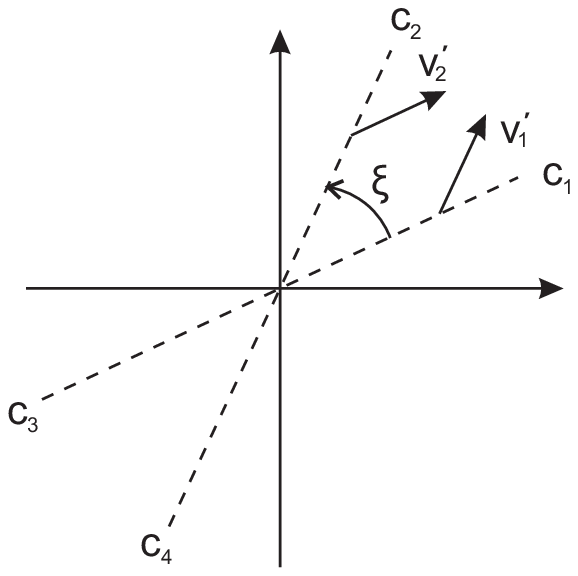}\\
  \caption{}
\end{figure}

 Put $\xi_1=\xi_3=\xi,\ \xi_2=\xi_4=\pi-\xi$, $\gamma_1=-\gamma_3=\varkappa_1,\ \gamma_2=-\gamma_4=-\varkappa_2$,  $v_1=(a/c,b/c)$, $v_2=(-b/c,-a/c)$, $ v_3=(-a/c,-b/c)$, $v_4=(b/c,a/c), $ $\theta_1=\theta_3=\xi-\pi/2,$  $\theta_2=\theta_4=\pi/2-\xi$. It is easy to check that  $(v_i,n_i)=1, i=1,2,3,4,$ where $n_i$ is the unit normal vector to $c_i$  that points anticlockwise. Indeed,
\begin{multline*}
(v_i,n_i)=\frac{\sqrt{a^2+b^2}}{c}\cos{(\pi/2-\xi)}=\frac{2}{2\sqrt{1-\alpha^2}}\sqrt{1-\alpha^2}=1,\\ i=1,2,3,4.
\end{multline*}
Equations (\ref{tx_1 prime}),(\ref{tx_2prime}) can be rewritten as follows
\begin{equation}
\tilde x(t)=\tilde x(0)+\tilde w(t)+\sum_{i=1}^4\gamma_iv_i L_{\tilde x}^{c_i}(t), \
t\leq\zeta,\label{eq4}
\end{equation}
where $\tilde x(t)=(\tilde x_1(t),\tilde x_2(t)),\tilde w(t)=(\tilde w_1(t),\tilde w_2(t)), \
\zeta$ is the first hitting time of 0 by the process
$(\tilde x(t))_{t\geq 0}$. So (\ref{eq4}) coincides with (\ref{main}).

Now we calculate the expression from Proposition for  $n=4$:
$$
S=\sum_{k=1}^4 \gamma_k\pi_k
\frac{\xi_{k-1}\xi_k}{{(\xi_{k-1}+\xi_k)+\gamma_k(\xi_{k-1}-\xi_k)}}\tan\theta_k.
$$

Let $\{\varkappa_1,\varkappa_2\}\in(-1,1)\setminus \{0\}$.  The invariant distribution $(\pi_i)_{i=1}^4$ can be obtained directly. But we make use of formula (29) from \cite{APi}. We have
$$
\pi_1=\frac{\tilde q_2\tilde q_3\tilde q_4}{\tilde p_1\tilde q_3\tilde q_4+\tilde p_1\tilde p_2\tilde q_4+\tilde p_1\tilde p_2\tilde p_3+\tilde q_2\tilde q_3\tilde q_4}.
$$
Taking into account (\ref{tp_k}),(\ref{tq_k}) we get
\begin{eqnarray*}
\pi_1 &=& \frac{(1-\gamma_2)\left((1+\gamma_1)(\pi-\xi)+(1-\gamma_1)\xi\right)}{D},\\
\pi_2 &=& \frac{(1+\gamma_1)\left((1-\gamma_2)(\pi-\xi)+(1+\gamma_2)\xi\right)}{D},\\
\pi_3 &=& \frac{(1+\gamma_2)\left((1-\gamma_1)(\pi-\xi)+(1+\gamma_1)\xi\right)}{D},\\
\pi_4 &=& \frac{(1-\gamma_1)\left((1+\gamma_2)(\pi-\xi)+(1-\gamma_2)\xi\right)}{D},\\
\end{eqnarray*}
where $D=2\left[(1+\gamma_1\gamma_2)\xi+(1-\gamma_1\gamma_2)(\pi-\xi)\right]>0.$ Then
$$
S=2\frac{\xi(\pi-\xi)}{D}(-\varkappa_1\varkappa_2\cot\xi).
$$
The condition $\xi\in(0,\pi)$ yields
$\cot\xi=-\frac{\alpha}{\sqrt{1-\alpha^2}}.$ Consequently $S>0$ if and
only if $\varkappa_1\varkappa_2\alpha>0.$

It is obvious that the processes $(x_1(t))_{t\geq 0}$ and
$(x_2(t))_{t\geq 0}$ meet in zero when and only when
$(\tilde x(t))_{t\geq 0}$ hits zero.

If  $\varkappa_1=\varkappa_2=1$ then there exists a unique invariant distribution
 $\pi_1=\pi_2=1/2, \ \pi_3=\pi_4=0$. It is easily to see that now
$S>0$ if and only if $\alpha>0.$

Finally, let $\varkappa_1=1, \ \varkappa_2\in(-1,1)\setminus\{0\}.$ Then the
invariant distribution is of the form $\pi_1=\tilde p_2/2,\ \pi_2=1/2,\
\pi_3=\tilde q_2/2,\ \pi_4=0.$ We get that $S>0$ if and only if
$\varkappa_2\alpha>0.$

The other cases when the modulus of at least one permeability coefficient is equal to 1 can be treated analogously.

Now let $\varkappa_1=0$. Then the unique invariant distribution is as follows $\pi_1=\pi_3=0,$ $\pi_2=\pi_4=1/2.$ It is easy to see that in this case $S=0$. Analogously, $S=0$ if $\varkappa_2=0$.

Thus we have proved the following statement.

\begin{Theorem}
Let $(x_1(0), x_2(0))\in\Bbb R^2\setminus\{(0,0)\}, \ \{\varkappa_1,\varkappa_2\}\subset[-1,1]$, $\alpha\in (-1,1)$. Then
\begin{eqnarray*}
1) \ \mathds P\{\exists\ t_0<\infty: x_1(t_0)=x_2(t_0)=0\}=1\ \hbox{if }\  \varkappa_1\varkappa_2\alpha>0,\\
2) \ \mathds P\{\exists\ t_0<\infty: x_1(t_0)=x_2(t_0)=0\}=0\
\hbox{if }\  \varkappa_1\varkappa_2\alpha\leq 0.
\end{eqnarray*}
\end{Theorem}

\begin{Remark} The conditions of processes meeting obtained in Theorem for $\alpha\in(-1,1)$ are completely different from those for $\alpha=1$ obtained in  \cite{B_local}.
\end{Remark}

\begin{Remark} As was mentioned above the process $(x(t))_{t\geq 0}$ almost surely does not hit any fixed point  except the points in which  at least one membrane is situated. It follows from statement 2) of Theorem that the process almost surely does not hit any fixed point in which exactly one membrane is situated.
\end{Remark}

\end{document}